\documentclass[12pt]{article}
\usepackage{amsmath,amsthm,amsfonts,amssymb,amscd}

\begin{document}

\theoremstyle{plain}
\newtheorem{thm}{\sc Theorem}[section]
\newtheorem{lem}[thm]{\sc Lemma}
\newtheorem{cor}{\sc Corollary}
\newtheorem{prop}[thm]{\sc Proposition}

\theoremstyle{remark}
\newtheorem{Case}{\bf Case}[section]
\newtheorem{rem}[thm]{\sc Remark}

\theoremstyle{definition}
\newtheorem{Def}{\bf Definition}[section]
\newtheorem*{pf}{\sc Proof}

\theoremstyle{definition}
\newtheorem{Conj}{\bf Conjecture}[section]
\newcommand{\epm}{{\cal E}_M}
\newcommand{\mim}{M \in \cal{M}}
\newcommand{\N}{{\mathbb N}}
\newcommand{\ve}{{\varepsilon}}
\newcommand{\be}{{\beta}}
\newcommand{\la}{{\lambda}}
\newcommand{\teta}{{\theta}}
\newcommand{\La}{{\Lambda}}
\newcommand{\ga}{{\gamma}}
\newcommand{\Ga}{{\Gamma}}
\newcommand{\po}{{\partial}}
\newcommand{\ov}{\overline}
\newcommand{\re}{{\mathbb{R}}}
\newcommand{\Pc}{{\mathbb{P}}}
\newcommand{\dc}{{\mathcal D}}
\newcommand{\na}{{\mathbb N}}
\newcommand{\co}{{\mathbb C}}
\newcommand{\fc}{{\mathcal F}}
\newcommand{\al}{{\alpha}}
\newcommand{\lc}{{\mathcal L}}
\newcommand{\om}{{\omega}}
\newcommand{\Z}{{\mathbb Z}}
\newcommand{\rth}{{\mathbb{R}^3}}
\newcommand{\W}{{\mathcal W}}
\renewcommand{\H}{{\mathcal H}}
\newcommand{\cc}{{\mathcal C}}
\newcommand{\E}{{\mathcal E}}
\newcommand{\PP}{{\mathcal{P}}}
\newcommand{\M}{{\mathcal M}}
\newcommand{\nep}{n(E)}
\newcommand{\de}{D(E)}
\newcommand{\ar}{A(R)}
\newcommand{\gp}{\pi}
\newcommand{\wtil}{\widetilde}
\newcommand{\st}{\subset}
\newcommand{\mo}{M_1}
\newcommand{\mto}{M_2}
\newcommand{\cz}{C^0}
\newcommand{\bbc}{\mathbb{C}}
\newcommand{\D}{\Delta}

\begin{title}{The topology, geometry and conformal structure of
    properly embedded minimal surfaces}\end{title}
\begin{author}
{Pascal Collin \and 
Robert Kusner \thanks{This research was supported by NSF grant DMS  0076085}
\and William H. Meeks, III \thanks{This research
 was supported by NSF grant DMS  0104044}
\and Harold Rosenberg} \end{author}
\maketitle
\section{Introduction}

Let $\M$ denote the set
of connected properly embedded minimal
surfaces in $\rth$ with at least two ends.
At the beginning of the past
 decade, there were two outstanding conjectures on the asymptotic geometry
 of the ends of an $M \in \M$ that were known to lead to topological
 restrictions on $M$.  The first of these conjectures, the 
 generalized Nitsche conjecture, stated that an annular end of such a $M \in
 \M$ is asymptotic to a plane or to the end of a catenoid.  Based on earlier
 work in \cite{mr6}, P.~Collin
\cite{col1} proved the generalized Nitsche conjecture.  In the case $\mim$
 has finite  topology,  the solution of this conjecture implies
that $M$ has finite total Gaussian curvature, which by previous work in  \cite{co3}, \cite{jm1}, \cite{lor1}, \cite{sc1} led to topological obstructions for such a minimal surface $M$. 

Our paper deals with the case where $M \in \M$ has infinite topology.
Before stating the second conjecture, we recall some definitions.  For
any connected manifold $M$, an end of $M$ is an equivalence class of
proper arcs on $M$ where two such arcs are equivalent if for any
compact domain $D$ in $M$, the ends of these arcs are contained in the
same noncompact component of $M-D$.  The set $\epm$ of all the ends of
$M$ has a natural topology that makes $\epm$ into a compact Hausdorff
space. The limit points in $\epm$ are by definition the {\itshape
limit ends} of $M$; an end $e\in\epm$ which is not a limit end will be
called a {\itshape simple end}.  To every $\mim$ is associated a
unique plane passing through the origin in $\rth$ called the {\itshape
limit tangent plane at infinity} of $M$ (see \cite{chm3}). The
existence of such a limit plane at infinity depends strongly on the
property that $M$ has at least two ends. For convenience, we will
always assume that the limit tangent plane at infinity is horizontal
or, equivalently, is the $x_1 x_2$--plane $P$.

A result of Frohman and Meeks \cite{fme2} states that the ends of $M$
are linearly ordered by their relative heights over $P$.  Furthermore,
they prove that this linear ordering, up to reversing it, depends only
on the proper ambient isotopy class of $M$ in $\rth$.  Since the space
of ends $\epm$ is compact and the ordering is linear, for any $\mim$,
there exists a unique {\itshape top end} which is the highest end in
the ordering on $\epm$. Similarly, the {\itshape bottom end} of $M$
is defined to be the end of $M$ which is lowest in the associated
ordering. The ends of $M$ that are neither top nor bottom ends are
called {\itshape middle ends} of $M$.

The second conjecture, motivated by analogy with the finite topology
setting, asserted that the middle ends of an $\mim$ are simple ends
which are $C^0$-asymptotic to a plane or to an end of a catenoid.
This conjecture was verified \cite{hm8} for middle ends of finite
genus, but remained open in the case of an infinite genus middle end,
where it was further conjectured that the limit to a plane or a
catenoid end must have finite integer multiplicity greater than one.
One consequence of this second conjecture is that the middle ends of
$\mim$ can be represented by proper subdomains with compact boundary
whose area in the ball $B_R$ of radius $R$ centered at the origin is
approximately equal to $n\pi R^2$ for some integer $n$ when $R$ is
large.  (Recall that a proper subdomain $E\subset M$ with compact
boundary is said to represent an end $e\in\epm$ if $E$ contains a
proper arc representing $e$.)

In this paper we will develop new fundamental theoretical tools for
understanding the topology, asymptotic geometry and conformal
structures of examples in $\M$. These tools are powerful enough to
prove that the middle ends of an $\mim$ are simple ends and have
quadratic area growth $n\pi R^2$.  An important consequence of these
methods is that the topology of examples in $\M$ with an infinite
number of ends is very restrictive.

\begin{thm}
If $\mim$, then a limit end of $M$ is a top or bottom end.  Thus $M$
has at most two limit ends, and in particular, $M$ can have only a
countable number of ends.
\end{thm}

Note that the above theorem gives strong topological restrictions that
a properly embedded minimal surface with an infinite number of ends
must satisfy.  (For example, the plane with a Cantor set removed has
an uncountable number of limit ends and so cannot properly minimally
embed in $\rth$.)  It is a consequence of the following geometric
result on the middle ends:

\begin{thm} \label{middle}
Suppose $\mim$. For a middle end $e$ of $M$, there is an associated
positive integer multiplicity $n(e)$. The multiplicity $n(e)$ is
defined by choosing a proper subdomain $E(e)$ with compact boundary
that represents $e$ such that the area $\ar = \mbox{Area}(B_R \cap
E(e))$ divided by $\gp R^{2}$ converges to $n(e)$ as $R \to
\infty$. Furthermore, $E(e)$ can be chosen so that for any other
representative $\wtil{E(e)} \st E(e)$ of $e\in\epm$ the associated
area function $\wtil{\ar}$ divided by $\gp R^{2}$ also converges to
$n(e)$ as $R \to \infty$.  In particular, a limit end of $\mim$ must
be a top or a bottom end.
\end{thm}

\noindent
Classical examples by Riemann \cite{ri1} and more recent examples by
Callahan, Hoffman and Meeks \cite{chm2} demonstrate that there exist
many $1$--periodic examples in $\M$ with two limit ends.

Theorem \ref{middle} is a crucial initial ingredient in the complete
topological classification theorem \cite{fme4} for properly embedded
minimal surfaces in $\rth$. Specifically, $\mo, \, \mto \in \M$ differ
by a diffeomorphism of $\rth$ if and only if they have the same genus
and, up to reversing the order of the ends, the corresponding ends
have the same genus (either $0$ or $\infty$) and the corresponding
integer multiplicities of the middle ends given in the statement of
Theorem \ref{middle} are the same modulo $2$.  Theorem \ref{middle}
has also played an essential role in the recent classification of
properly embedded nonsimply connected periodic minimal surfaces of
genus zero \cite{mpr1}.

We will also show that some of the examples in $\M$ have strong
restrictions on their conformal structure as well. Recall that a
Riemannian surface $M$ is {\itshape recurrent} if almost all Brownian
paths are dense in $M$.  An important conformal property for recurrent
Riemannian surfaces is that positive harmonic functions are constant.

\begin{thm} \label{recur}
If $\mim$ and $M$ has two limit ends, then $M$
is recurrent.
\end{thm}

Since triply-periodic minimal surfaces have one end and are never
 recurrent, some restriction on the number of ends is necessary for
 the conclusion of Theorem \ref{recur} to hold. Embeddedness is also a
 necessary hypothesis in Theorem \ref{recur}, since at the end of
 Section 3 we will construct properly immersed minimal surfaces with
 two ends in $\rth$ that have nonconstant bounded harmonic functions.
 Indeed, it has been conjectured that if $\mim$, then $M$ is recurrent
 \cite{me23}.

The proof of Theorem \ref{recur} depends on a basic result (Theorem
\ref{thm3.1}) on the conformal structure of a properly immersed minimal
surface with boundary contained in a closed halfspace of $\rth$.
Theorem \ref{thm3.1} has played an important role in the proof of
uniqueness of the helicoid \cite{mr8} and also
in the proof of the invariance of flux for a coordinate function of a
properly immersed minimal surface (see \cite{me23}).

Our paper is organized as follows. In Section~2 we apply the classical
Weierstrass representation of minimal surfaces to derive some special
proper superharmonic functions defined on certain regions of a
properly immersed minimal surface. Next, we use these special
functions, together with the divergence theorem, to prove that middle
ends of an $\mim$ have quadratic area growth. In Section 3 we again
use these special functions to derive some of our basic theorems on
conformal structure.

\section{Quadratic area growth of middle ends}

In the proof of the ordering theorem \cite{fme2}, Frohman and Meeks 
prove that every middle end of a surface $\mim$ is contained between
two catenoids in the sense of the following definition.

\begin{Def} 
Suppose $M$ is a properly immersed minimal surface with compact boundary in
 $\rth$.  We will say that $M$ is contained between two catenoids if for some $c_1>0$, $M\subset\{(x_1,x_2,x_3)
 \mid |x_3|\le c_1 \ln r, \, r^2=x_1^2+x_2^2, \,  r \geq 2 \}$. \end{Def}

Since the middle ends of a properly embedded minimal
  surface are contained between two catenoids, the following lemma implies
  that the middle ends of a properly embedded minimal surface
  are never limit ends.

\begin{lem}\label{quad} If $M$ is a properly immersed minimal surface
 with compact boundary and $M$ is contained between two catenoids,
 then $M$ has quadratic area growth. This means that the area of $M$
 in the ball $B_R = \{(x_1,x_2,x_3)\mid x_1^2+~x_2^2+x_3^2\le R^2\}$ is at  most $C\pi(R + 1)^2$
 for some positive $C$. Furthermore, such an $M$ has at most $C$ ends.
\end{lem}

The proof of the above lemma depends on a fundamental inequality given in
 the next lemma. (In \cite{me24} the calculations in Lemma
 \ref{lem3.2} were repeated to obtain related results for ends of periodic minimal surfaces.)

\begin{lem}\label{lem3.2}
Let $M$ be a minimal surface and assume $r=\sqrt{x_1^2 + x_2^2}\ne 0$ on
 $M$. Then, $|\Delta_M\ln r|\le |\nabla_M x_3|^2/r^2$.
\end{lem}

\begin{pf}
Assume $M$ is not a plane and note that the points where the gradient
 of $x_3$ is zero 
 are isolated on $M$; so it suffices to prove the inequality stated in the
 lemma on the complement of the horizontal points.

Let $g$ denote the stereographic projection of the Gauss map of $M$ to the
 extended complex plane $\bbc \cup \infty$.  Then, by the classical
 Weierstrass representation, the coordinates $(x_1,x_2,x_3)(z)$ are given by
$$
\text{Re}\int \left(\frac12(\frac1g-g),\quad \frac{i}{2}(\frac1g+g),1
 \right)\, dz,
$$
where $z=x_3+ix_3^*$.  We will let $\Delta_z$ denote the planar
Laplacian in $z$-coordinates.  

Then
$$
(x_1+ix_2)(z)=\ov{\int\frac{1}{2g}}dz-\int \frac{g}{2}dz=\ov \xi-\mu,
$$
where
$
\xi'=\frac1{2g},\quad\mu'=\frac{g}{2},\quad \xi'\mu'=\frac14.
$
Letting $\ln w = \ln r + i \theta$ for $x_1 + ix_2 = w =
re^i{^\theta}$, we have
$$
\Delta_z \ln(x_1 + ix_2) = 4\frac \partial{\partial\ov z} \frac
\partial{\partial z}\ln(\ov
\xi-\mu)=4\frac \partial{\partial\ov z}
 \left( \frac{-\mu'}{\ov\xi-\mu} \right)
=
 \frac{4\mu'\ov\xi'}{(\ov \xi-\mu)^2} $$
$$= \Delta_z\, (\ln |\ov \xi-\mu|+i \,\Delta_z \theta) 
= \Delta_z\, \ln \, r+i\, \Delta_z \theta.$$
Thus,
$
\Delta_z\, \ln\, r= \text{Re}\left (4\frac{\mu'\ov\xi'}{(\ov\xi-\mu)^2}
 \right)= \text{Re} \left( \frac{g}{\ov{g}}\cdot \frac1{r^2 e^{2 i \theta}}
 \right).$  Since $|\frac {g}{\ov{g}}\cdot \frac1{e^{2 i \theta}}| = 1,$ then \,
$|\Delta_z\, \ln \, r|\le \frac1{r^2}$.  Since $\Delta_M = |\nabla_M
x_3|^2 \Delta_z,$ this  completes the proof of Lemma 2.2. 
\qed
\end{pf}

\noindent
{\bf Proof of Lemma 2.1:} Let $ C_t=\{p\in \rth \mid r(p)=t\}$ be the
vertical cylinder of radius $t$ and let $M_t$ be the part of $M$
inside $C_t$.  Since the part of $M$ inside the ball of radius $R$
centered at the origin is contained in $M_R$, it suffices to prove
that $M_R$ has quadratic area growth as a function of $R$.

In the complement of the $x_3$-axis,   one has the ordered orthonormal basis:
 \\
 $(\nabla r,\nabla
 x_3,r\nabla\theta)=(A_1,A_2,A_3)$. Let $B_i=A_i-(\vec n\cdot A_i)\vec n$
 be the tangent part of $A_i$ (here $\vec n$ is the unit normal to
 $M$), so
$(\vec{n} \cdot A_1)^2 + (\vec{n} \cdot A_2)^2 + (\vec{n} \cdot A_3)^2
 = 1$ and
$$
|B_i|^2=|A_i|^2-(\vec n\cdot A_i)^2=1-(\vec n\cdot A_i)^2.
$$
Hence $|B_1|^2+|B_2|^2=1 + (\vec n \cdot A_3)^2 \geq 1$. Since $B_1 =
\nabla_M r$ and $B_2 = \nabla_M x_3$, 
$$
|\nabla_Mr|^2+|\nabla_Mx_3|^2\ge1.
$$
Thus,
$$
\int_{M_R}dA\le\int_{M_R}(|\nabla_Mr|^2+|\nabla_M x_3|^2)\, dA.
$$
Therefore, it remains to prove that both $\int_{M_R}|\nabla_Mr|^2\, dA$  and
 $\int_{M_R}|\nabla_M x_3|^2\, dA$ grow at most quadratically in $R$.

Without loss of generality, after removing a compact subset  and
 homothetically scaling $M$, we may assume that the third coordinate
 function on $M$ satisfies the inequality $|x_3| \leq \frac12 \ln r.$ 
 Consider the function $f\colon M \to \re$ defined by
$-x_3\arctan(x_3) + \frac12 \ln(x_3^2 + 1)$. A calculation yields
$\D_M (f) = \frac{-|\nabla_M x_3|^2}{x_3^2 +1}.$

By Lemma \ref{lem3.2}, $\Delta_M \ln r \leq |\nabla_M x_3|^2/r^2$ and so the
 function $h=\ln r + f(x_3)$ is superharmonic on $M$.  Since $\ln r$
 is proper in a closed region of $\mathbb{R}^3$ containing $M$ and $M$
 is proper, it follows that $\ln r$ is a proper function on $M$. And
 $h \geq \frac{1}{10} \ln r$, so $h$ is a proper nonnegative superharmonic function on $M$.

Now for any positive proper $C^2$-function $H$ on $M$ and $T\geq \sup
(H(\po M)),$
$$
\int_{H^-{^1}{([0,T])}}\Delta_M H= -\int_{\po M} \nabla_M H \cdot \eta
 +\int_{H^-{^1}(T)}| \nabla_M H |,
$$ 
where $\eta$ is the outward pointing conormal to the boundary.  Hence,
 if $\Delta_M H\le0$, then $\int_{H^-{^1}(T)}|\nabla_M H|$ is positive
 monotonically decreasing as $T\to\infty$. So $\Delta_M H\in L^1(M)$
 and, choosing $H=h$, we see that $\D_M h \in L^1(M).$ Since for $r$
 large, $|\D_M h| \geq \frac12| \D_M f|$, we note that $\D_M f \in
 L^1(M)$ as well.  Hence,
$$
\int_{M_R}|\Delta_M f| \, dA = \int_{M_R} \frac{|\nabla_M x_3|^2}{x_3^2 +
 1}\, dA \leq c_2
$$
for some constant $c_2$.  But $|x_3| \leq \ln(R)$ on $M_R$, so we have
$$
\int_{M_R} |\nabla_M x_3|^2 \, dA \leq \int_{M_R} 
\left( \frac{(\ln R)^2+1} {x_3^2 + 1} \right) |\nabla_M x_3|^2 \, dA
\leq [(\ln R)^2+1]c_2. 
$$
This completes the proof that $\int_{M_R} |\nabla_M x_3|^2 \, dA $ grows at
 most quadratically in $R$.

Since $\D_M f \in L^1(M)$ and $|\D_M f| \ge |\Delta_M\ln r|$, we also have
$\Delta_M \ln r\, \in \, L^1(M)$.
Because
$$
\int_{M_R} \D_M \ln r =-\int_{\po M} \frac{\nabla_M r \cdot \eta}{r} +
 \int_{C_R \cap M} \frac{|\nabla_M r|}{R} = c_3 + \frac{1}{R}\int_{C_R\cap
   M} |\nabla_M
 r|
$$
 and $\int_{M_R} |\D_M \ln r|$ converges, $\frac1R \int_{C_R \cap M} |\nabla_M r|$
 has a finite limit as $R \to \infty.$  Hence 
$$\int_{C_R \cap M} |\nabla_M r| \le c_4 R
$$
 for some constant $c_4$.
In fact
$\int_{C_{\rho} \cap M}|\nabla_M r| \le c_4 \rho$
for any $\rho \in [1, R]$, so
 the coarea formula implies
$$
\int_{M_R}|\nabla_M r|^2\le \int_1^R c_4 \rho d\rho \leq \frac{c_4 R^2}2 \quad 
$$
 which means that $\int_{M_R}|\nabla_M r|^2$ grows quadratically in $R$.  It 
follows that the area of $M$ grows
 quadratically in $R$.

We now check that $M$ has a finite number of ends.  If $E$ is a proper
noncompact subdomain of $M$ with $\partial E$ compact, then we know
that the area of $E$ grows at most quadratically.  But, by the 
monotonicity formula for area \cite{si1}, the area of $E$ must grow
asymptotically at least as quickly as the area of a plane, which means
that for large $R$ the area of $E$ inside the ball $B_R$ is at least
$\pi(R-1)^2$, and if $M$ has at least $n$ ends, then, for every
$\varepsilon > 0$,  the area of $M_R$ must be greater than 
$(n-\varepsilon)\pi(R+1)^2$ for large $R$.  In particular, the last sentence
in the statement of Lemma \ref{quad} holds, which completes our
proof.  \qed

\vspace{.17in}
\begin{rem}In the statement of Lemma \ref{quad}, the hypothesis that
  $M$ lies between two catenoids can be weakened to the property that
  $M$ lies above a catenoid end and intersects some positive vertical
  cone in a compact set.  Under this weaker hypothesis, the conclusion
  of Lemma \ref{quad} that $M$ has quadratic area growth still holds.
  To prove this more general result one uses the function $h_c = h
  +cx_3$ for some positive $c$ in place of the function $h$ defined in
  the proof of the lemma.  Indeed, under this weaker hypothesis $h_c$
  is again proper and positive for sufficiently large $c$; since $M$
  is minimal, $\D_M x_3 = 0$ and so $|\D_M h_c|=|\D_M h|\geq \frac12|
  \D_M f|$, implying $\D_M f \in L^1(M)$ as before.  We also have
  $|x_3| \leq aR$ on $M_R$ for some positive $a$, and so, with 
  constant $c_2$ chosen as above,
$$
\int_{M_R} |\nabla_M x_3|^2 \, dA \leq \int_{M_R} 
\left( \frac{(aR)^2+1} {x_3^2 + 1} \right) |\nabla_M x_3|^2 \, dA 
\leq [(aR)^2+1]c_2.
$$
The rest of the argument is unchanged.

\end{rem}


We now explain how Theorem \ref{middle} stated in the Introduction
follows from Lemma \ref{quad}. Let $M$ be as in the statement of Lemma
\ref{quad} and let $P$ denote the $x_1 x_2$-plane.  The monotonicity
formula for area \cite{si1} implies that
$\lim_r{_\rightarrow}{_\infty} A(r)/\pi r^2$ exists and is a finite
number $n(M)$.  Since $M$ can be viewed as a locally finite integral
varifold with compact boundary and $M$ has quadratic area growth,
standard compactness theorems (see \cite{si1}) imply that the sequence
of integral varifolds $\frac{1}{k} M = \{(\frac{x_1}{k},
\frac{x_2}{k}, \frac{x_3}{k}) \mid (x_1, x_2, x_3) \in M \}$ converges
to the locally finite integral varifold $n(M) P$ as $k \rightarrow
\infty$.  Hence, the area-multiplicity $n(M)$ of $M$ is a positive
integer.  The integer $n(M)$ can be easily identified from the proof
of Lemma \ref{quad} as $n(M) = \lim_R{_\rightarrow}{_\infty}
\frac{1}{2\pi R} \int \limits _{M \cap C_R} |\nabla_M r|$.  (This
identification will be used later in the proof of Theorem \ref{thm3.3}
and Lemma \ref{lem3.1}.)  Then Theorem \ref{middle} stated in the
Introduction follows immediately from these comments and Lemma
\ref{quad}.

\section{Parabolicity}

We now apply the results of the previous section to derive some global
 results on the conformal structure of properly immersed or properly
 embedded minimal surfaces. We will say that a Riemannian surface $M$
 with boundary is {\it parabolic\/} if bounded harmonic functions on
 the surface are determined by their boundary values. We recall that
 given a point $p$ on a Riemannian surface $M$ with boundary, then
 there is an associated measure $\mu_p$ on $\po M$, called the
 ``hitting'' or {\it harmonic\/} measure, such that $\mu_p(I)$, for an
 interval $I\subset\po M$, is the probability that a Brownian path
 beginning at $p$ ``hits'' the boundary a first time at a point in
 $I$.  Note that harmonic measure enjoys a domain monotonicity
 property: if $M'\subset M$ is a subdomain containing $p$, and if also
 $I\subset\po M'$, then the corresponding harmonic measures satisify
 $\mu'_p(I)\leq\mu_p(I)$; this is because the family of Brownian paths
 from $p$ to $I$ within $M'$ is contained in the corresponding family
 of paths within $M$.

It is well known that $M$ is parabolic if and only if the harmonic
 measure $\mu_p$ for any $p\in \text{Int}(M)$ is full, that is,
 $\int_{\po M}d\mu_p=1$. In fact if $\mu_p$ is full and $f\colon
 M\to\re$ is a bounded harmonic function, then for any $p\in
 \text{Int}(M)$, $f(p)=\int_{\po M}f(x)d\mu_p$.  It is easy to check
 that if $\mu_p$ is full for some point $p\in \text{Int}(M)$, then
 $\mu_q$ is full for any other point $q\in \text{Int}(M)$.

In order to verify whether a Riemannian surface $M$ with boundary is
 parabolic it is sufficient to find a proper nonnegative superharmonic
 function $h\colon M\to[0,\infty)$.  To see this suppose
 $f_1,f_2\colon M\to\re$ are two bounded harmonic functions on $M$
 with the same boundary values and $f_1(p)>f_2(p)$ for some
 $p\in\text{Int}(M)$. Then consider the proper function $H_t\colon
 M\to\re$ defined by $H_t(x)=h(x)-t(f_1(x)-f_2(x))$. For $t$
 sufficient large, $H_t(p)<0$ and hence $H_t$ has a minimum at some
 interior point of $M$, contradicting the minimum principle for
 superharmonic functions.

With this preliminary discussion in mind we now state the first theorem of
 this section.

\begin{thm}\label{thm3.1}
If $M$ is a connected properly immersed minimal surface  in $\re^3$, possibly with
 boundary, then $M(+)=\{(x_1,x_2,x_3)\in M\mid x_3\ge0\}$ is parabolic.
\end{thm}

\begin{pf}
Let $M(n) = \{(x_1, x_2, x_3) \in M \mid 0 \leq x_3 \leq n \}$.
We first prove that
$$
M(n,*)=\{(x_1,x_2,x_3)\in M(n) \mid 1 \le x_1^2+x_2^2\}
$$
is parabolic. Let $r(x_1,x_2,x_3)=\sqrt{x_1^2+x_2^2}$ and define $h\colon
 M(n,*)\to[0,\infty)$ by $h(p)=\ln(r(p))-x_3^2(p)$. By Lemma \ref{lem3.2}, $h$ is a
 superharmonic function on $M(n,*)$ and is proper since $\ln r$ is proper
 and $x_3^2$ is bounded. Because $h$ is eventually positive, $M(n,*)$ is parabolic. As $M(n)$ is the
 union of $M(n,*)$ and the compact surface $M(n)\cap\{r\le 1\}$,
 $M(n)$ is also parabolic.

We now check that $M(+)$ is parabolic by proving that each component
 $C$ of $M(+)$ is parabolic. Let $p\in C$ be a point with positive
 third coordinate; by rescaling we will assume $x_3(p)=1$.  Since
 $M(n)$ is parabolic, so is $C(n) = C \cap M(n) = C\cap
 x_3^{-1}([0,n])$, and thus the relation $1 = x_3(p)$ can be evaluated
 as an integral:
$$
1 = \int_{\po(C(n))}x_3(x) d\mu_p(n) 
$$
$$
  = 0\cdot\int_{\po C\cap x_3^{-1}(0)} d\mu_p(n) +
    \int_{\po C(n)\cap x_3^{-1}((0,n))} x_3(x) d\mu_p(n) + 
    n\cdot\int_{\po C(n)\cap x_3^{-1}(n)} d\mu_p(n),
$$
where $\mu_p(n)$ is the harmonic measure on the boundary of $C(n)$.
The middle term is nonnegative, so 
$$
\int_{\po C(n)\cap x_3^{-1}(n)} d\mu_p(n)\le\frac1n.
$$ 
Since $\mu_p(n)$ is full on $C(n)$,
$$
\int_{\po C(n)-x_3^{-1}(n)} d\mu_p(n)\ge 1-\frac 1n.
$$
\noindent
Taking limits (using domain monotoncity of $\mu_p(n)$) as $n \rightarrow
\infty$, one obtains $\int_{\po C} d\mu_p=1$, which proves the
theorem.
\end{pf} \qed


\begin{rem}Recently Meeks \cite{me23} has applied Theorem \ref{thm3.1}
  to prove that the flux of a coordinate function of a properly
  immersed minimal surface is well defined.  Also, Meeks and Rosenberg \cite{mr7} have used this
  Theorem \ref{thm3.1} to prove that if $M$ is a finite topology
  properly immersed minimal surface in $\rth$ such that a plane
  intersects $M$ transversely in a finite number of component curves,
  then $M$ is a conformally a finitely punctured Riemann surface. This
  result and Theorem \ref{thm3.1}  should be compared with the theorem of
  Morales \cite{moral1} which proves the existence of a proper
  conformal minimal immersion of the open unit disk into $\rth$.
\end{rem}

\begin{cor}\label{cor3.1}
Suppose $D$ is a proper domain in $\re^2$ and $M$ is a minimal graph over
 $D$ which is bounded from below. Then $M$ is parabolic.

\end{cor}

Recall that a complete Riemannian surface $M$ is called \it{recurrent}
 \rm for
 Brownian motion if, with probability one, a Brownian path starting at a
 point $p\in M$ will enter every neighborhood of any other point $q\in M$
 for a divergent sequence of times. The notion of being recurrent is closely
 related to the notion of parabolicity for surfaces with boundary. If $M$ is the union of two
subdomains that intersect in a compact subset of $M$, then $M$ is recurrent
 if these subdomains are parabolic. Since a properly immersed minimal
 surface $M$ in $\re^3$ can be expressed as $M=M(+)\cup M(-)$, where
 $M(-)=\{(x_1,x_2,x_3)\in M \mid x_3\le 0\}$, and $M(+)$ and $M(-)$ are
 parabolic, then if $M(+)\cap M(-)=M\cap x_3^{-1}(0)$ is compact, $M$ is
 recurrent. We restate this result as a corollary.

\begin{cor}\label{cor3.2}
If $M$ is a properly immersed minimal surface in $\re^3$ and some plane
 intersects $M$ in a compact set, then $M$ is recurrent for Brownian motion.
\end{cor}

In certain cases it can be shown that a properly embedded minimal surface
 has a compact intersection with some plane. The final theorem of this
 section gives an important instance of this compact intersection property.

\begin{thm}\label{thm3.3}
If $M$ is a properly embedded minimal surface with two limit ends, then
 between any two middle ends of $M$, there is a plane that intersects $M$
 transversely in a compact set. In other words, given two distinct middle
 ends of $M$ there is a plane $P$ that intersects $M$ transversely in a
 compact set and the representatives of these ends in $M-P$ lie on opposite
 sides of $P$. In particular, $M$ is recurrent for Brownian motion.
\end{thm}

The above theorem will follow immediately from the next lemma.

\begin{lem}\label{lem3.1}
Suppose that $M$ is a noncompact properly immersed minimal surface with
 compact boundary contained between vertical catenoid ends $C_1$ and
 $C_\la$, where $C_1$ has logarithmic growth 1 and $C_\la$ has logarithmic
 growth $\la\ge1$ and $C_\la$ lies above $C_1$. If the asymptotic
 area growth of $M$ is $n\pi r^2$, then the vertical flux $F= \int_{\po
 M} |\nabla_M x_3|$ satisfies $\quad 2\pi n\le F \le 2\pi n\la$.
\end{lem}

\begin{pf} 
Recall the function $f = -x_3 \arctan (x_3) + \frac{1}{2} \ln
(x_3^2+1)$ defined in the proof of Lemma \ref{quad}.  In that proof it
was shown that $\D_M f, \, \D_M \ln \,  r \in L^1(M)$.  It also follows
from that proof for every real number $c$, $h_c= \ln \, r + cx_3 + f$ is
superharmonic outside of a compact subdomain of $M$.  Under the
hypotheses given in the statement of Lemma \ref{lem3.1}, for a
divergent sequence of points $\{p_i\}$ in $M$,
$\lim_i{_\rightarrow}{_\infty} \,  h_c (p_i) = + \infty$ for $ c>
\frac{\pi}{2}-{\frac{1}{\lambda}}$ and $\lim_i{_\rightarrow}{_\infty} \,  h_c(p_i) = -
\infty$ for $c< \frac{\pi}{2} - 1$.

Assume now that $c > \frac{\pi}{2}-{\frac{1}{\lambda}}$. The type of calculations carried out in the proof of Lemma
\ref{quad} imply that for 
$T$ sufficiently large, $\int \limits_{h_c^-{^1}(T)} \nabla_M h_c
\cdot \eta$ is a positive monotonically decreasing function of $T$,
where $\eta$ is the outward pointing unit conormal to
$h_c{^-}{^1}((-\infty, T])$.
(In the case $c < \frac{\pi}{2} -1$, $\int \limits_{h_c^-{^1}(T)}
\nabla_M h_c \cdot \eta$ is a negative and monotonically increasing  in
norm as a function of $-T$).
Hence, for $c > \frac{\pi}{2} -\frac{1}{\lambda}$ and $T$ large, \,  
$$0 \leq \int \limits _{h_c^-{^1}(T)} \nabla_M h_c \cdot \eta = \int \limits
_{h_c^-{^1}(T)} \frac{\nabla_M r}{r} \cdot \eta + \int \limits
_{h_c^-{^1}(T)} \nabla_M (cx_3 + f) \cdot \eta.$$  
Since $\nabla_M f = -\nabla_M x_3 \, \cdot \arctan (x_3)$ and $x_3 \rightarrow
\infty$ as $T \rightarrow \infty$, 
$$
(c-\frac {\pi}{2}) \int \limits
_{\partial M} \nabla_M x_3 \cdot \eta = (\frac {\pi}{2}-c) \lim_{T \rightarrow \infty} \int \limits
_{h_c^-{^1}(T)} \nabla_M x_3 \cdot \eta \leq \lim_{T \rightarrow \infty} \int \limits
_{h_c^-{^1}(T)} \frac {\nabla_M r} {r} \cdot \eta .
$$
Let $c=\frac{\pi}{2}-\frac{1}{\lambda}$ and we obtain $$
F = \int
\limits _{\partial M}|\nabla_M x_3| = - \lim_{T \rightarrow \infty} \int \limits
_{h_c^-{^1}(T)} \nabla_M x_3 \cdot \eta \leq \lambda \cdot \lim_{T \rightarrow \infty} \int \limits
_{h_c^-{^1}(T)} \frac {\nabla_M r} {r} \cdot \eta.
$$
But since $\D_M \ln \, r \in L^1 (M),$ the divergence theorem
implies
 $$\lim _{T
\rightarrow \infty}  \int \limits _{H_c^-{^1}(T)} \frac {\nabla_M r}
{r} \cdot \eta = \lim _{R \rightarrow \infty} \int \limits _{M \cap
  C_R} \frac {\nabla_M r}{R} \cdot \eta = \lim _{R \rightarrow \infty}
\frac {1}{R} \int \limits _{M \cap C_R} | \nabla_M r |,$$ 
 where $C_t$ is the cylinder of radius $t$ centered along the
 $x_3$-axis.  From the discussion immediately following the proof of Lemma \ref{quad},
$$
\lim _{R \rightarrow \infty} \frac{1}{R} \int \limits _{M \cap C_R} |
\nabla_M r | = 2 \pi n .
$$ 
Hence, for $c = \frac {\pi}{2} -\frac{1}{\lambda}$, we obtain $F
\leq 2 \pi n \lambda.$   Making similar
calculations in the case $c < \frac {\pi}{2} - 1,$ we obtain the
inequality $F\geq 2 \pi n$, which completes the proof of the lemma.

\end{pf} \qed

\noindent {\bf Proof of Theorem 3.3:}

\rm  Suppose $M$ has two limit ends with horizontal limit tangent
plane at infinity. From the proof of the ordering theorem \cite{fme2}, there exists
an end $E$ of a vertical catenoid or of a horizontal plane between any two middle
ends of $M$.  If $E$ were a vertical
catenoid with positive logarithmic growth between middle ends $e_1$
and $e_2$, where $e_2$ is the next middle end of $M$ above $e_1$, then
between $e_2$ and the end $e_3$ just above $e_2$, there would be an
end $\tilde{E}$ of a catenoid between $e_2$ and $e_3$ such that the
logarithmic growth of $\tilde{E}$ is at least equal to the logarithmic
growth of $E$.

After a homothety of $\rth$, we may assume that the logarithmic growth
of $E$ is $1$.  By Lemma \ref{lem3.1}, the flux of $\nabla_M x_3$ across
the boundary of any proper domain $M(e_2)$ representing $e_2$ must be at
least $2\pi$.  Similarly, for the end $e_3$, the flux of $\nabla_M x_3$
across the boundary of any proper domain $M(e_3)$ representing $e_3$
must also be at least $2\pi$.  In fact, if we let $\{e_2, e_3, \hdots,
e_n, \hdots\} \subset \epm$ be the infinite set of middle ends above $e_1$, and
ordered so that $e_i < e_j$ if $i<j$, then for any such end $e_j$
there exist proper disjoint subdomains $M(e_i)$ representing the $e_i$ such that the
flux of $\nabla_M x_3$ across $\partial M(e_i)$ is at least $2 \pi$. 

Assume now that $P$ is a horizontal plane that intersects $E_1 =
M(e_1)$ in a circle $S^1 \subset P$ which we may assume is the
boundary of $E_1$.  Let $D$ be the disk in $P$ with $\partial D =
S^1$.  Without loss of generality, we may assume that $D$ intersects
$M$ transversely in a finite number of simple closed curves which separate
$M$ into a finite number of components.  Let $M^\prime$ be one of the
components above $D \cup E_1$ that contains an infinite number of
middle ends of $M$.  Then $x_3 \colon M^\prime \rightarrow [x_3(P),
\infty)$ is proper and so the flux of $\nabla_M x_3$ across $\partial
M^\prime$ is at least as big as the sum of the fluxes of $\nabla_M x_3$
coming from the middle ends of $M^\prime$.  Since each middle end of $M^\prime$
contributes at least $2 \pi$ of flux, the total flux of $\nabla_M x_3$
across $\partial M^\prime$ must be infinite.  But $\partial M^\prime$
is compact and the flux is no more than the total length of the
boundary of $M^\prime$ which is finite.  This contradiction proves $E$
must be a horizontal plane from which
Theorem \ref{thm3.3} follows. \qed

\vspace{.17in}

\noindent
{\bf An example with non-constant bounded harmonic functions:}

It is conjectured that every $\mim$ is recurrent, and that, for
a properly embedded minimal surface $M$ with one end, every positive
harmonic function on $M$ is constant.\footnote{About 20 years ago,
Dennis Sullivan asked whether a positive harmonic function on a
properly embedded minimal surface in $\rth$ must be constant.}  It is
known that a positive harmonic function on any properly embedded
triply-periodic minimal surface $f \colon N \rightarrow \rth$ is
constant, but such an $N$ is never recurrent.  We will now construct a
two sheeted covering space with $p \colon M \rightarrow N$ which has
non-constant harmonic functions, and so $f\circ p \colon M \rightarrow
\rth$ is a properly immersed minimal surface with two ends and with
non-constant bounded harmonic functions.

First recall that $N$ is topologically an infinite genus surface with
one end.  Let $\Gamma = \{a_1, a_2, \hdots, a_n, \hdots \}$ be a countable
proper collection of closed curves on $N$ which generate the first
homology group of $N$.  Let $\sigma \colon \pi_1(N) \rightarrow
\mathbb{Z}_2$ be the homomorphism which factors through
$\tilde{\sigma} \colon H_1(N) \rightarrow \mathbb{Z}_2$ where
$\tilde{\sigma}([a_1]) = 1$ and $\tilde{\sigma} ([a_i]) = 0$ for $i
\neq 1$.  Let $p\colon M \rightarrow N$ be the $\mathbb{Z}_2$-cover of
$N$ corresponding to Ker$(\sigma) \subset \pi_1(N)$.  Let $D_1 \subset D_2 \subset \hdots$ be a proper
compact exhaustion of $N$ by subdomains with one boundary curve and
such that $a_1 \subset D_1$.  Let $E = N-$Int$(D_1)$ and $E(+)$ and
$E(-)$ be the two components of $p^-{^1}(E)$.  Note that $E(+)$ and
$E(-)$ are proper domains which represent the two ends of $M$.  Let
$\tilde{D_1} \subset \tilde{D_2}\subset \hdots$ be the associated
pullback compact exhaustion of $M$ and note that each domain
$\tilde{D_i}$ has two boundary curves, $\partial (i, -) \subset E(-)$
and $\partial (i, +) \subset E(+)$, respectively.

Let $h_n \colon \tilde{D}(n) \rightarrow [-1, 1]$ be the harmonic
function with boundary value $-1$ on $\partial (n, -)$ and $+1$ on
$\partial (n, +)$.  Since $\{h_n\}$ is a uniformly bounded, increasing
sequence, it converges to a harmonic function $h\colon M
\rightarrow [-1, 1]$.  We will now prove that $h$ is nonconstant by
showing $h$ has $-1$ as an asymptotic limiting value on the end of
$E(-)$ and has $+1$ as an asymptotic limiting value on the end of
$E(+)$.  

Consider a divergent sequence $p(i) \in N$ such that $p(i) \in D(i) -
D(i-1)$.  Let $W(n) = D(n) - $Int$(D_1)$.  For $n>i$ consider $P =
\int \limits _{\partial D_1} \mu (p(i), n)$ where $\mu(p(i),n)$ is the
hitting measure for $p(i)$ considered to lie in $W(n)$.  Since $N$ is
not recurrent for Brownian motion, for every $\varepsilon >0$, there
exists an $N(\varepsilon)$ such that if $n>i>N(\varepsilon)$, then the
probability $P$ of a Brownian path starting at $p(i)$ in $W(n)$ of
exiting a first time at $\partial D_1$, is less than $\varepsilon$.
This implies that if $p(i,+)\in E(+)$ and $p(i,-) \in E(-)$ are the
two lifts of $p(i)$ to $M$, then $h_n(p(i,+)) \geq 1-2 \varepsilon$
and $h_n(p(i,-)\leq -1 +2 \varepsilon$, which implies $h(p(i), +) \geq
1-2 \varepsilon$ and $h(p(i,-1) \leq -1 + 2 \varepsilon$.  By letting
$\varepsilon \rightarrow 0$, we obtain our earlier claim that $h$ is
asymptotic to $-1$ on the end of $E(-)$ and asymptotic to $+1$ on the
end of $E(+)$.

\bibliographystyle{plain}
\bibliography{bill}

\end{document}